\documentclass{nmeauth}

\usepackage{bm}
\usepackage{amsmath}
\usepackage{amsfonts}
\usepackage{graphicx}
\usepackage{rotating}
\usepackage{float}
\usepackage[ruled,vlined]{algorithm2e}
\usepackage{caption}
\usepackage{color}
\definecolor{blue}{rgb}{0.,0.,1.}
\definecolor{red}{rgb}{1.,0.,0.}

\newcommand{\figref}[1]{Figure \ref{#1}}

\usepackage{hyperref}
\hypersetup{
    colorlinks=true,
    linkcolor=blue,
    filecolor=magenta,      
    urlcolor=blue,
    citecolor=blue,
    pdftitle={Local Element Operations for Curved Simplex Meshes},
    }
    
\begin{document}
\NME{1}{20}{XX}{40}{23}

\runningheads{A. Shi and P.-O Persson}{Local Element Operations for Curved Simplex Meshes}

\title{Local Element Operations for Curved Simplex Meshes}

\author{A. Shi\affil{1},\affil{2}\comma\corrauth, P.-O. Persson\affil{1}}
\corraddr{A. Shi, Department of Mathematics, University of California, Berkeley, Berkeley, CA 94720-3840, U.S.A. (andrewshi@math.berkeley.edu)}
\address{\affilnum{1}\ Department of Mathematics, University of California, Berkeley, CA, USA \\
\affilnum{2}\ Department of Mathematics, New York University Shanghai, Shanghai, China}

\noreceived{...}
\norevised{}
\noaccepted{}

\begin{abstract}
Mesh optimization procedures are generally a combination 
of node smoothing and discrete operations 
which affect a small number of elements
to improve the quality of the overall mesh.
These procedures are useful as a post-processing step 
in mesh generation procedures 
and in applications such as fluid simulations
with severely deforming domains. 
In order to perform high-order mesh optimization, 
these ingredients must also be extended to high-order (curved) meshes.
In this work, we present a method to perform 
local element operations on curved meshes.
The mesh operations discussed in this work are edge/face swaps, 
edge collapses, and edge splitting (more generally refinement) 
for triangular and tetrahedral meshes.
These local operations are performed by first identifying 
the patch of elements which contain the edge/face being acted on, 
performing the operation as a ``straight-sided one" by
placing the high-order nodes via an isoparametric
mapping from the master element, 
and smoothing the high-order nodes on the elements 
in the patch by minimizing a Jacobian-based high-order mesh distortion measure.
Since the initial ``straight-sided guess" from the placement 
of the nodes via the isoparametric mapping frequently results 
in invalid elements, the distortion measure must be regularized 
which allows for mesh untangling for the optimization to succeed.
We present several examples in 2D and 3D to demonstrate 
these local operations and how they can be combined with 
a high-order node smoothing procedure to maintain mesh quality 
when faced with severe deformations. 
\end{abstract}

\keywords{high-order; curved mesh; mesh optimization; edge flip; edge collapse; deforming domains}

\section{INTRODUCTION}
%%Why high order and curved meshing.
High-order methods have received considerable interest in the 
computational science community in the last two decades. 
This is due to their potential to deliver highly accurate solutions 
at a comparable cost to traditional low-order methods. 
In fact, for certain classes of problems --- such as turbulent flow --- 
it is widely believed that high-order methods will be \textit{required} 
for accurate and grid-converged solutions \cite{deville2002high}.
Much of the great promise of high-order methods is due their ability to 
handle complex geometries with unstructured meshes.
These methods require high-order (curved) meshes to ensure the curved geometry 
is represented with sufficient accuracy to enable 
the high-order convergence rate of these methods.
However, the existence (or rather, lack thereof) of 
robust and automated high-order mesh generation procedures 
is still one of the main bottlenecks in the widespread adoption 
of these methods in the design process.
To this end, the high-order community has identified high-order mesh generation 
as a topic deserving of considerable research attention 
and as one of the top pacing items in high-order method research \cite{wang2013high, wang2016perspective}.
One might say curved meshing first distinguished 
itself as a distinct research topic within mesh generation
in 2015, the first year it had its own dedicated session 
at the International Meshing Roundtable \cite{imr24}. 
In terms of available software, Gmsh is a widely used free software 
with the capability to generate and display curved meshes 
for elements of arbitrarily high-order \cite{geuzaine2009gmsh}. 
Up until fairly recently, there were no commercial high-order mesh generators 
available before Pointwise released support for this in 2019 \cite{pointwise}. 
We refer the reader to the introduction of \cite{turner2017high} 
for a recent survey of the current state of 
affairs of high-order meshing research and additional references.

%%Mesh cleanup phases.
Most straight sided mesh generation procedures proceed in two steps: 
1) the initial creation of a mesh from geometry through 
a method like advancing front, Delaunay or octree, and 
2) a mesh optimization post-processing step to improve element quality. 
There is a large body of literature on these 
mesh quality improvement procedures. 
While those works may differ in their 
notions of optimality and specific techniques used, 
they largely agree mesh optimization via node movement
and discrete, localized element operations 
are the best way to perform mesh quality improvement.
Freitag and Ollivier-Gooch \cite{freitag1997tetrahedral} 
consider various node smoothing procedures 
and face and edge swapping operations, 
and show that each mechanism fails to significantly improve 
the mesh quality when used individually 
but result in very high quality meshes when combined. 
They offer a variety of empirical recommendations on how to do so 
and demonstrate the effectiveness of their schedules 
on a variety of tetrahedral mesh geometries.
De Cougny and Shephard \cite{de1999parallel} combine edge collapses 
and splitting operations for the purpose of mesh adaptation. 
Klinger and Shewchuk \cite{klingner2008aggressive} 
consider a richer set of operations 
in an effort to take these ideas to an extreme 
to find the highest quality mesh, 
assuming that speed is not the highest priority. 
We refer the reader to Section 3 of their paper 
where they survey many of the commonly 
(as well as lesser used) discrete mesh operations.
These are some of the ``classic" papers in the field for 
mesh quality improvement for linear meshes, 
and while the focus of this paper is on curved meshes, 
the field of linear mesh quality improvement is still highly active.
Some more recent contributions include, but are certainly not limited to 
directions like new discrete operations to reach better local minima
\cite{dassi2016tetrahedral, chen2016shell}, and the inclusion of parallelization
\cite{freitag1999parallel,sastry2014parallel,lopez2022parallel}.

%High order smoothing - current
Now the question is how to extend these ideas to curved meshes. 
Existing work \cite{remacle2014optimizing, karman2016high, dobrev2019target, gargallo2015optimization} 
on high-order node smoothing primarily takes the approach of 
optimization-based minimization of some energy functional or objective function 
that quantifies high-order mesh distortion. 
The key issues here are the need for mesh distortion metrics 
for high-order elements \cite{roca2011defining, gargallo2014defining} 
and a procedure for high-order untangling \cite{gargallo2015optimization, toulorge2013robust, stees2020untangling}.
In principle, some of the approaches to curved mesh generation 
such as the solution of a linear/nonlinear elasticity analogy 
\cite{persson2009curved, xie2013generation, nielsen2002recent, moxey2014thermo} 
or a PDE-based approach \cite{fortunato15winslow} 
could also be adapted for the purpose of high-order mesh smoothing.
There is practically no literature on the use of local element operations for curved meshes.
To the best of the authors' knowledge, the only existing work on 
high-order mesh operations \cite{feuillet2018p} considers 
face and edge swaps for tetrahedral meshes with second-order ($P^2$) elements.

It is clear that further developments in these areas will be required 
to advance the current state of curved mesh generation and adaptivity. 
The difficulty of directly acting on high-order meshes is one possible explanation for why nearly all curved mesh generation
approaches are of the \textit{a priori} variety
rather than direct methods \cite{shontz2021direct} which proceed in the same 
two phases as most straight-sided mesh generators do.
In the context of finite element solutions, solution quality is 
frequently dictated by the mesh element of the worst quality.
Local mesh operations for curved meshes could be an inexpensive and direct way 
to deal with this and other applications such as the repair of 
deforming high-order domains occurring in time-dependent fluid simulations.

\begin{figure}[ht!]
    \centering
    \includegraphics[width=0.75\textwidth, angle=0]{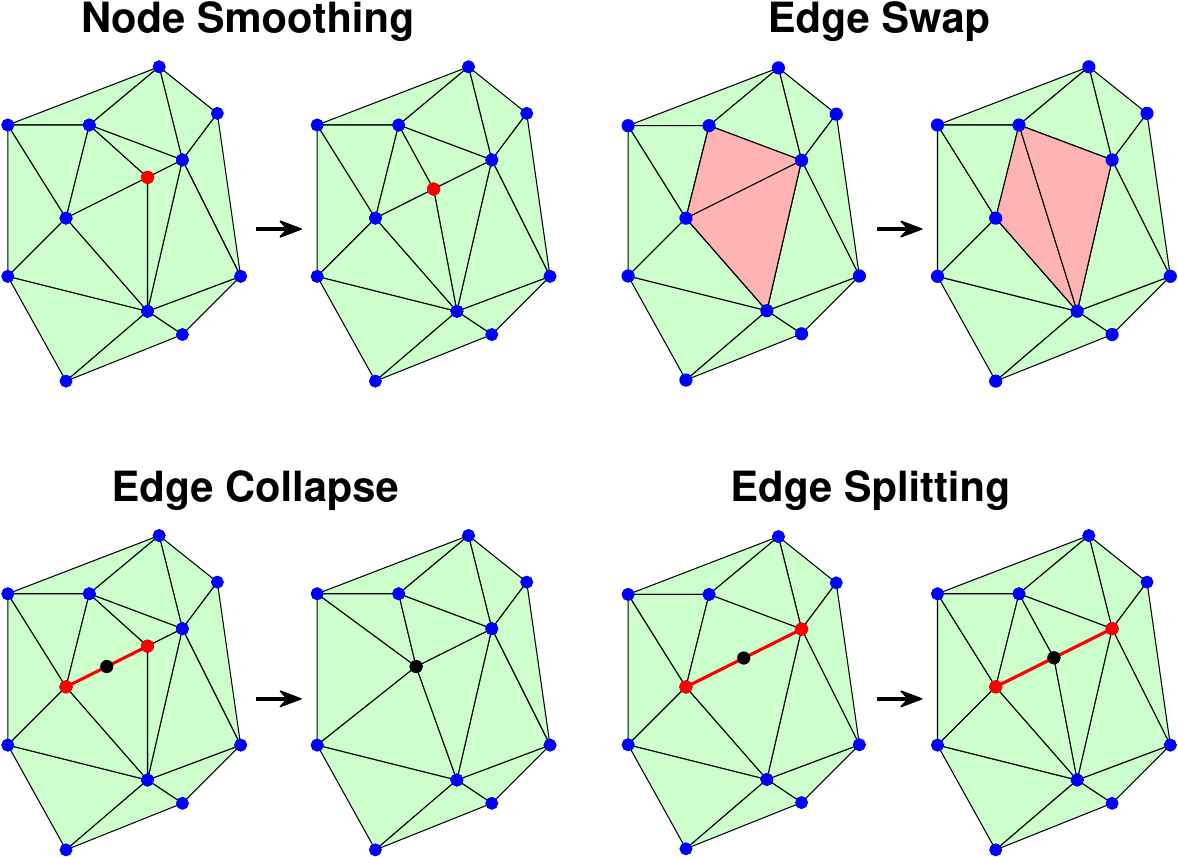}
    \caption{Straight sided mesh operations.}
    \label{fig:p1ops}
\end{figure}
\begin{figure}[ht!]
    \centering
    \includegraphics[width=0.75\textwidth, angle=0]{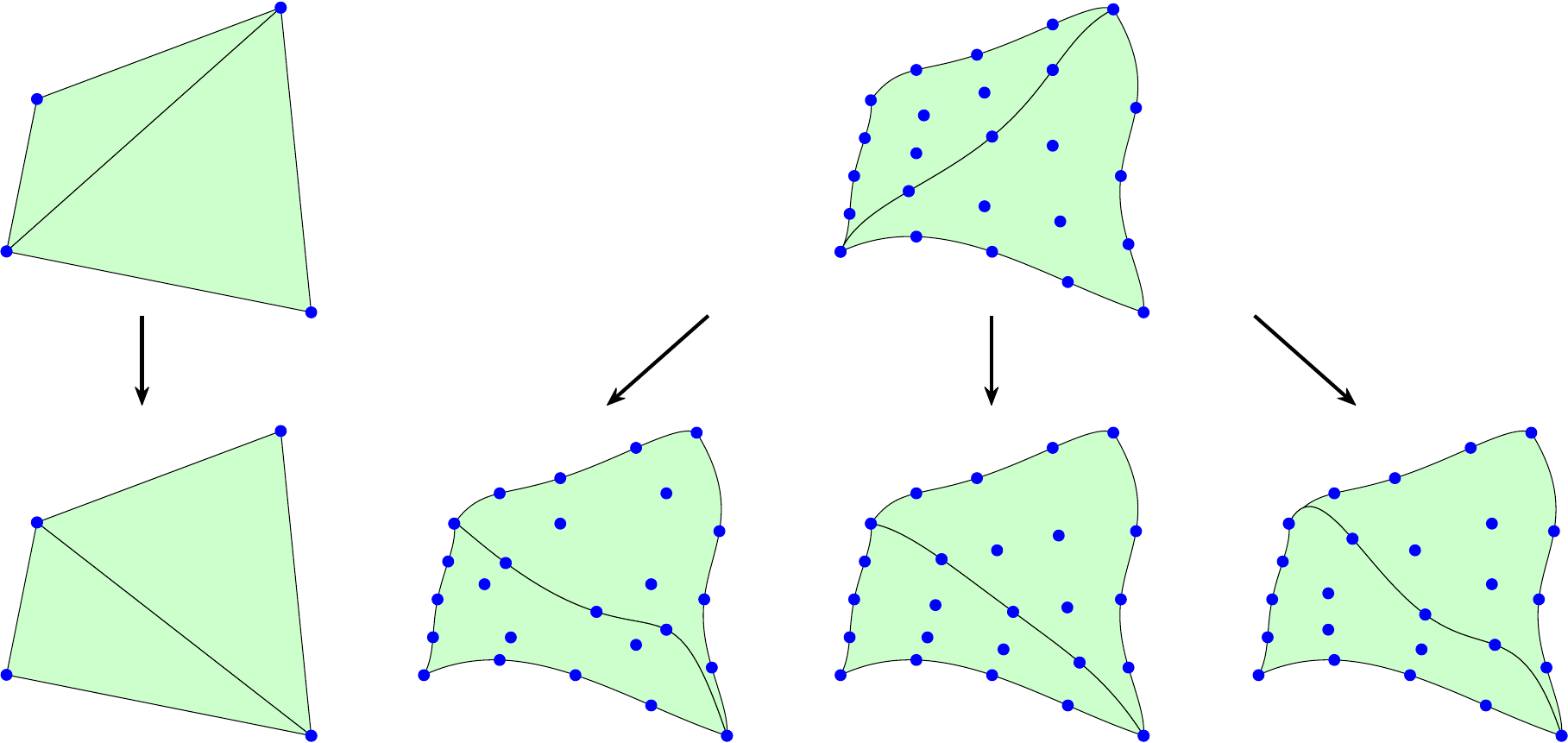}
    \caption{Only one possible flip for straight-sided elements \textit{(left)}, but many possible flips for high-order elements \textit{(right)}.}
    \label{fig:manyhighpflips}
\end{figure}

%Main goal of paper and structure
The main goal of this paper is to answer the question, 
``How we can extend local element operations for
straight-sided meshes (\figref{fig:p1ops}) 
to high-order triangular and tetrahedral meshes?"
One key difference from the low-order case is that there are many valid ways 
to perform these operations due to the additional degrees of freedom 
afforded by the high-order nodes (\figref{fig:manyhighpflips}). 
So there is now an \textit{optimal} way to perform these local operations 
by placing these additional nodes in a way to minimize 
a high-order Jacobian-based mesh distortion measure.
In particular, this paper will discuss high-order face and edge swaps, 
edge collapses and edge splitting for triangular and tetrahedral meshes.
The remainder of this article is organized as follows. 
In Section 2, we introduce the high-order Jacobian based distortion measure 
and its untangling capabilities. 
In Section 3, we describe the process of performing these operations, 
and discuss some practical details related to their implementation.
In Section 4, we present multiple examples illustrating how these operations 
can be combined with high-order node smoothing 
to repair curved meshes that undergo severe deformations.
Finally, in Section 5, we offer concluding remarks and directions for future research.

\section{DISTORTION MEASURE AND UNTANGLING}

We describe the construction of the high-order distortion
measure used in this work for node smoothing 
from a standard linear distortion measure,
which is then regularized to 
allow for simultaneous smoothing and untangling \cite{branets2002distortion,escobar2003simultaneous}.
We also offer some comments on how this measure
is used in our setting for optimization.

\subsection{High-order distortion measure}
%a) Background on the linear measure
We introduce the distortion measure for linear elements (\ref{EQ: DISTORTION MEASURE}). 
This measure -- known as the inverse mean-ratio shape measure -- 
was first introduced in \cite{liu1994shape, liu1994relationship}.
This measure is widely used in the meshing community as a 
shape distortion metric due to its many desirable properties 
such as the standard invariance to translation, rotation, and scaling, 
but is also special as a geometric shape measure that also enjoys
many of the desirable properties of an \textit{algebraic} mesh quality metric \cite{knupp2001algebraic}. 
In particular, this measure is shown to be 
convex \cite{munson2007mesh} and well-suited for optimization. 

As an aside, we note the metric we are about to introduce is geometric in nature, 
as are most metrics in use today, as opposed to being PDE or solution based. 
In the course of using mesh generation and adaptation in the numerical solution of PDEs, 
we point out that using a metric that is geometric in nature may not perfectly align 
with the ultimate goal of solving the PDE since the metric employed was not based on the PDE solution. 
This is addressed in \cite{knupp2007remarks} where the author discusses the notion of 
``bridging the gap between a priori quality metrics and solution-dependent metrics". 

\begin{figure}[ht]
    \centering
    \includegraphics[width=0.4\textwidth,angle=0]{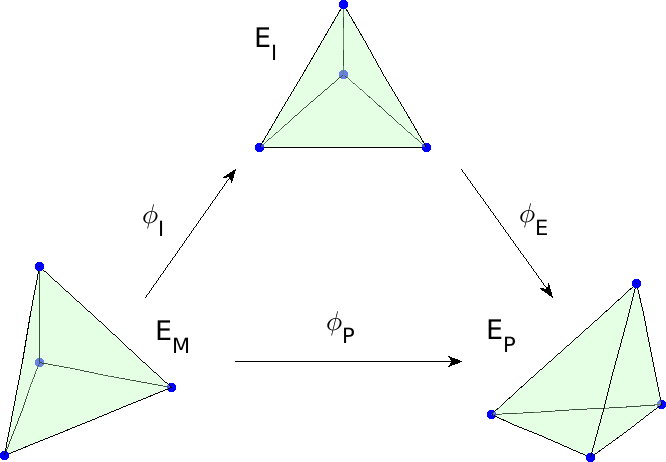}\hspace{1cm}\includegraphics[width=0.4\textwidth, angle=0]{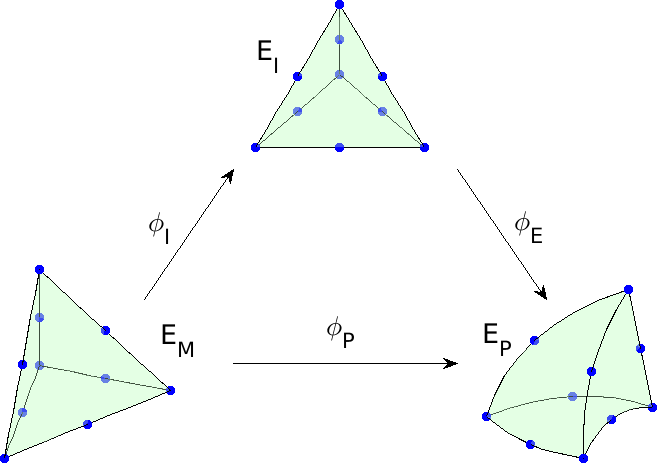} 
    \caption{Mappings between the master, ideal and physical elements: (a) linear case \textit{(left)} and (b) high-order case \textit{(right)}.}
    \label{fig:mappings}
\end{figure}

%b) Define the linear measure
Consider a linear element $E^I$ having the desired shape and size (the ideal element) 
and a corresponding element in physical space $E^P$. 
From these we can define a unique affine mapping $\phi_E: E^I \rightarrow E^P$. 
It is frequently convenient to introduce the notion of the master element $E^M$ 
and define two additional mappings $\phi_I, \phi_M$ (\figref{fig:mappings}a) that map from the master element to the ideal and physical element, respectively. 
We can define $\phi_E$ in terms of the these mappings as follows
\begin{equation} \label{EQ: PHI_E}
    \phi_E: E^I \xrightarrow{\phi_I^{-1}} E^M \xrightarrow{\phi_P} E^P.
\end{equation}
The element shape distortion metric $\eta$ is defined as
\begin{equation} \label{EQ: DISTORTION MEASURE}
    \eta(D\phi_E) = \frac{||D\phi_E||_F^2}{d|\sigma|^{2/d}}
\end{equation}
where $d$ is the spatial dimension, $||\cdot||_F$ is the Frobenius norm, 
and $\sigma = \det(D\phi_E)$. Note that for linear elements, 
since $\phi_E$ is affine, its Jacobian $D\phi_E$ is constant. 
This distortion measure quantifies the deviation of the shape 
of the physical element $E^P$ with respect to the ideal element $E^I$. 
This distortion measure is 1 for the ideal element and tends to infinity 
as the physical element degenerates. 
It is often convenient to define a corresponding quality measure $q$ as
\begin{equation} \label{EQ: QUALITY DEFINITION}
    q = \frac{1}{\eta} \in [0, 1].
\end{equation}

%c) Define the high-order measure
For high-order elements, we can no longer apply (\ref{EQ: DISTORTION MEASURE}) 
directly because the mapping $\phi_E$ is no longer affine, 
so $D\phi_E$ is no longer constant. 
Roca et al. describes an extension of the linear distortion measure 
to high order triangular \cite{roca2011defining} 
and tetrahedral \cite{gargallo2014defining} elements. 
The idea is that since the linear measure quantifies the
local deviation between the ideal and physical elements, 
we can obtain a high-order element-wise distortion $\hat{\eta}$ 
measure by integrating the linear measure on the physical high order element
\begin{equation} \label{EQ: HIGH ORDER DISTORTION}
    \hat{\eta}(D\phi_E) = \left(\frac{1}{|t|}\int_t \eta^r (D\phi_E)  \ dx \right)^{1/r}
\end{equation}
where $|t|$ is the area of the physical element and $r = 2$ in this work. 
Explicit expressions for the isoparametric mappings 
for the high-order case (\figref{fig:mappings}b) in terms of shape 
functions are given in \cite{roca2011defining}.

\subsection{Regularization and untangling}
%d) Untangling
The distortion measure $\eta$ needs to be modified
to eliminate the singularity so an optimization procedure 
can recover from an invalid initial configuration. 
We use the regularization procedure introduced independently
by \cite{branets2002distortion, escobar2003simultaneous} around the same time,
but follow the notation \cite{escobar2003simultaneous}
which replaces $\sigma$ in the 
denominator of (\ref{EQ: DISTORTION MEASURE}) by 
\begin{equation} \label{EQ: REG JACOBIAN}
    \sigma_{\delta}(\sigma) = \frac{1}{2}\left(\sigma + \sqrt{\sigma^2 + 4\delta^2} \right)
\end{equation}
where $\delta$ is a positive element-wise parameter. 
This regularized Jacobian $\sigma_{\delta}(\sigma)$ 
is a monotonically increasing function of $\sigma$ such that 
$\sigma_{\delta}(0) = \delta$, which tends to 0 when $\sigma$ 
tends to $-\infty$, which allows us to overcome the vertical asymptote 
at $\sigma = 0$ in the original measure. 
This parameter $\delta$ is only set for a non-zero value 
if there exists an invalid element in the mesh under consideration, 
otherwise it is set to zero for a mesh with all valid elements. 
Using the regularized Jacobian, we can now modify (\ref{EQ: DISTORTION MEASURE}) 
to obtain a linear distortion measure capable of 
\textit{simultaneous} smoothing and untangling
\begin{equation} \label{EQ: REG DISTORTION MEASURE}
    \eta_{\delta}(D\phi_E) = \frac{||D\phi_E||_F^2}{d|\sigma_{\delta}|^{2/d}}
\end{equation}
and use this to define a high-order regularized measure $\hat{\eta}_{\delta}$ 
analogous to (\ref{EQ: HIGH ORDER DISTORTION}). 
Figures 5 and 6 of \cite{gargallo2015optimization} provide a useful illustration 
of how this regularized distortion metric compares to the original for a simple example.

\begin{figure}[ht!]
    \centering
    \includegraphics[width=0.5\textwidth, angle=0]{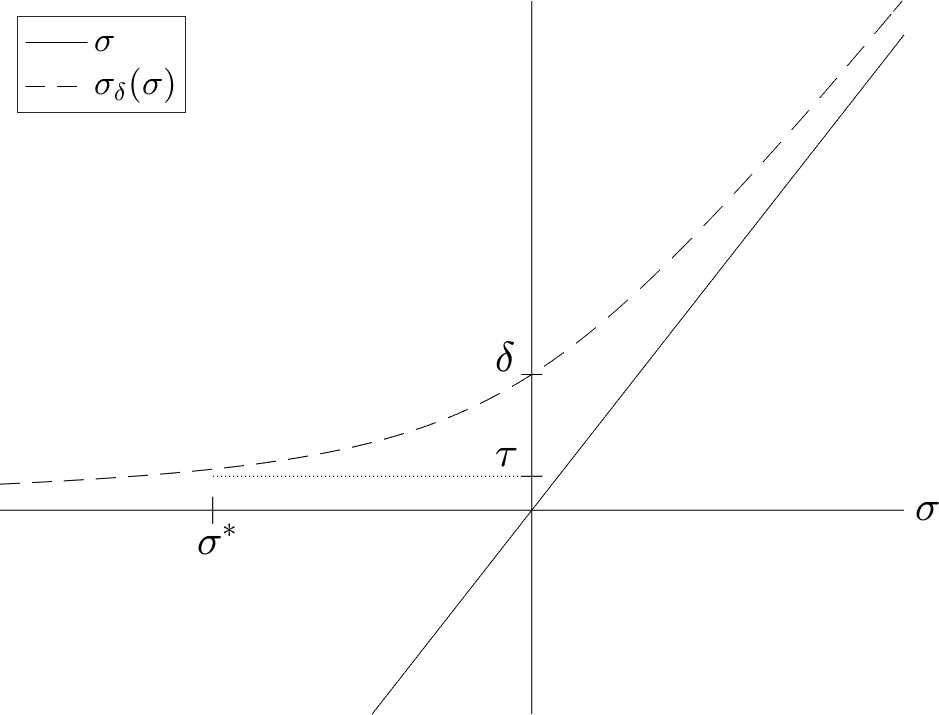} \\
    \caption{Representation of $\sigma_{\delta}(\sigma)$ (reproduced from \cite{gargallo2015optimization}).}
    \label{fig:regdist}
\end{figure}

%e) how to pick the parameter
What remains is the issue of how to pick this 
positive elementwise regularization parameter $\delta$. 
We certainly require the original ($\hat{\eta}$) and regularized ($\hat{\eta}_{\delta}$)
distortion measures to have nearby minima, 
so $\delta$ needs to be sufficiently small.
On the other hand, $\delta$ has to be large enough 
to ensure that $4\delta^2$ is significant compared to $\sigma^2$ 
in the expression for $\sigma_{\delta}$.
When this idea was originally introduced in \cite{escobar2003simultaneous}, 
this parameter was chosen in a fairly ad hoc manner 
by testing a few values for a given initial tangled mesh.
In \cite{gargallo2015optimization}, a heuristic was
developed to choose a constant value of $\delta$ for each element. 
Their setting deals with high-order mesh generation
by curving the boundaries of a well-shaped straight-sided mesh 
and minimizing the regularized distortion metric. 
For them, in the mapping $\phi_E$, 
each physical element $E^P$ in the curved mesh has 
the ideal element $E^I$ as the corresponding element
in the original straight-sided mesh.
They propose to choose this parameter $\delta$ 
solely based off information from the straight-sided ideal element.
Defining $\sigma^* = -\det \phi_I$ and imposing
\begin{equation}
\sigma_{\delta}(\sigma^*) = \frac{1}{2}\left(\sigma^* + \sqrt{(\sigma^*)^2 + 4\delta^2} \right) = \tau > 0
\end{equation}
for some given tolerance $\tau$ implies
\begin{equation}
\delta(\sigma^*) = \frac{1}{2}\sqrt{(2\tau + |\sigma^*|)^2 - (\sigma^*)^2} =\sqrt{\tau^2 + \tau |\sigma^*|}.
\end{equation}
They argue that $\tau$ should be small compared to $\sigma$ and select
$\tau = \alpha|\sigma^*|$, giving the final value for $\delta$ as
$$\delta(\sigma^*) = |\sigma^*|\sqrt{\alpha^2 + \alpha}.$$ 
The value $\alpha = 10^{-3}$ is observed to work well in practice
and accomplish the tradeoff required on the value of $\delta$
(Figure \ref{fig:regdist}).
Our setting deals with high-order mesh smoothing,
so we choose the ideal element $E^I$ to always be the 
equilateral triangle/tetrahedron.
For the unregularized distortion measure, the actual 
size of $E^I$ is irrelevant since the measure 
is invariant to scaling. This is no longer true 
for the regularized distortion measure. 
So in order to maintain the assumption from their setting that each 
$E^I$ and $E^P$ are roughly of the same size, 
we need to scale the ideal element to have the same
volume as each corresponding physical element and choose $\delta$
accordingly.

\subsection{Optimization}
\label{sec:optim}
%f) how to aggregate
We need to aggregate the distortion measures 
for individual elements into a single quantity 
for a group of elements for optimization purposes.
What we really want is to optimize for the worst element, 
since the distortion of the worst element in a mesh
is often what drives the solution quality in the context of numerical simulation.
But since minimizing the maximum distortion is difficult
optimization-wise, we consider the normalized sum of squares
of the elementwise distortion of each element $e_i$
in a group of elements $M$
\begin{equation}\label{EQ: AGG}
    \hat{\eta}_{agg}(M) = \frac{1}{|M|}\sum_{e_i \in M}\hat{\eta}_{\delta}(D\phi_{e_i})^2
\end{equation}
where $\hat{\eta}_{\delta}$ is the regularized high-order distortion measure 
obtained by combining (\ref{EQ: HIGH ORDER DISTORTION}) and (\ref{EQ: REG DISTORTION MEASURE})
\begin{equation} \label{EQ: REG HIGH ORDER DISTORTION}
    \hat{\eta}_{\delta}(D\phi_E) = \left(\frac{1}{|t|}\int_t \eta_{\delta}^r (D\phi_E)  \ dx \right)^{1/r}.
\end{equation}

\begin{figure}[h]
    \centering
    \includegraphics[width=0.3\textwidth]{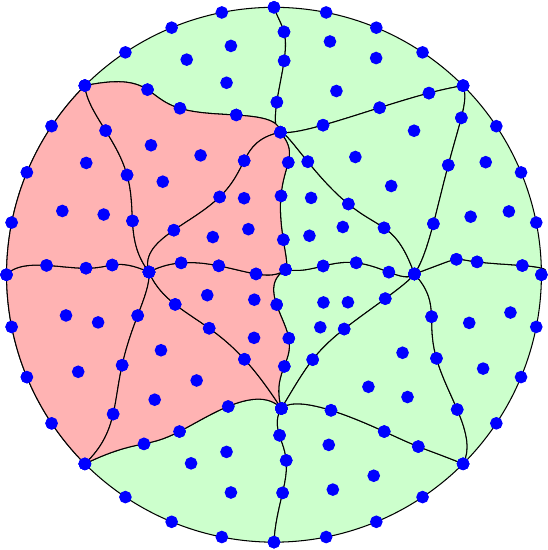}
    \hspace{2mm}
    \includegraphics[width=0.3\textwidth]{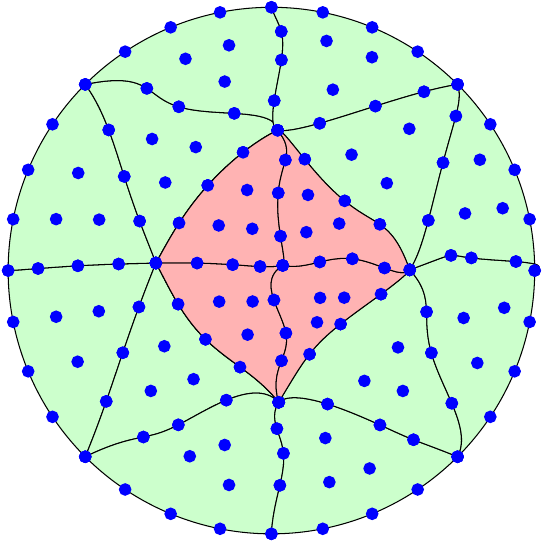}
    \hspace{2mm}
    \includegraphics[width=0.3\textwidth]{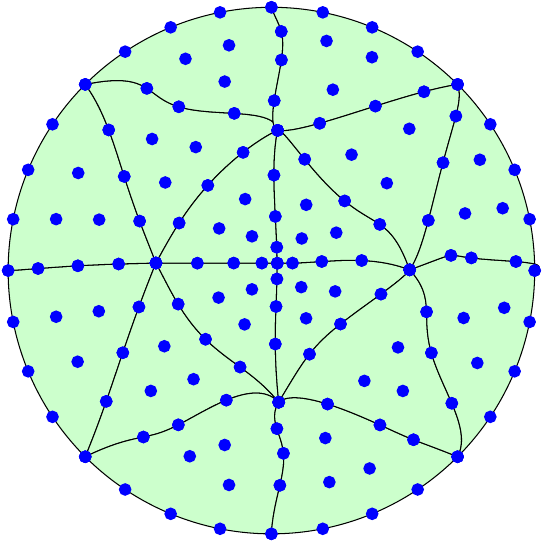}
    \caption{Patch-based high-order mesh smoothing.}
    \label{fig:mshpatch}
\end{figure}

%g) how to do smoothing
For smoothing over the entire mesh, we largely follow the
localized approach given in Appendix A of \cite{gargallo2015optimization}, 
which considers a ``Gauss-Seidel" like method which 
iteratively updates each node by a Newton step with
backtracking line search until either the objective function 
or the amount of node movement satisfies some stopping criteria.
Instead of looping over nodes, we loop over patches
(Figure \ref{fig:mshpatch}), which seems to help
mitigate some of the ``back and forth" the optimization
procedure can empirically observed to become stuck in, 
which is also observed in \cite{dobrev2019target} for a 
different high-order mesh smoothing procedure.
We note this solver offers no guarantees for especially
tangled or pathological meshes due to 
the issues with the regularization method which many authors
have pointed out fails quite often in practice 
such as if the boundary of the mesh is also tangled \cite{irving2004invertible},
or the possible failure of the high-order Gauss quadrature rule
to detect inversion at the corners of our high-order elements.

\section{ELEMENT OPERATIONS}
We describe how to perform the local element operations
of edge/face swaps, edge collapses, and edge splitting for
curved, simplex elements in 2D and 3D. 
The operations roughly proceed by freezing the boundary
of the patch of affected elements, performing the operation
as a ``straight-sided one", and 
then simultaneously untangling and smoothing the result.
We note some practical issues along with differences
in the high-order setting.

\subsection{Edge/face swaps}
\begin{figure}[h!]
    \centering
    \includegraphics[width=0.2\textwidth]{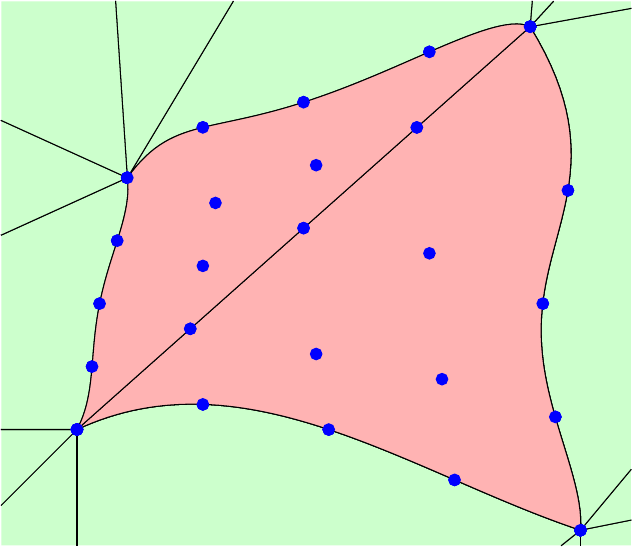} 
    \hspace{2mm}
    \includegraphics[width=0.2\textwidth]{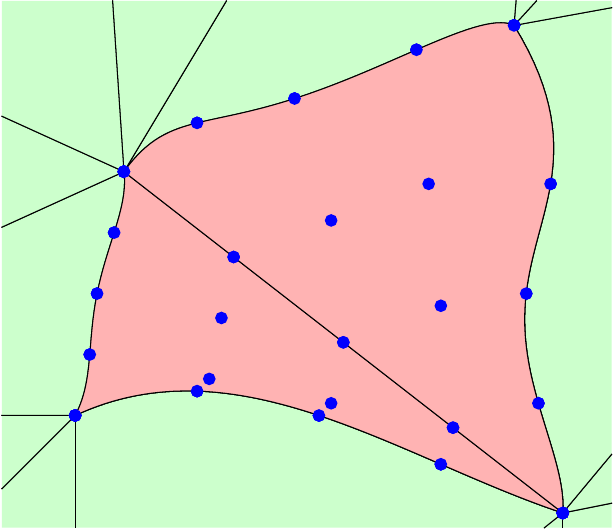} 
    \hspace{2mm}
    \includegraphics[width=0.2\textwidth]{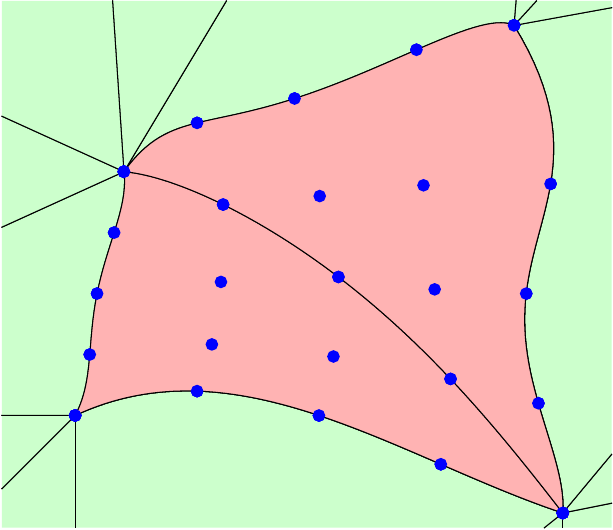} 
    \caption{Edge flips in 2D: original elements \textit{(left)}, initial straight-sided flip \textit{(middle)}, curved flip \textit{(right)}.}
    \label{fig:2dflip_schem}
    \includegraphics[width=0.23\textwidth]{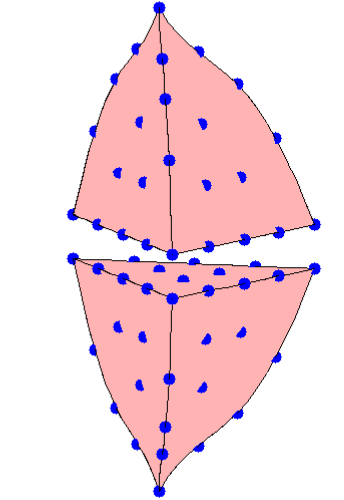} 
    \hspace{2mm}
    \includegraphics[width=0.23\textwidth]{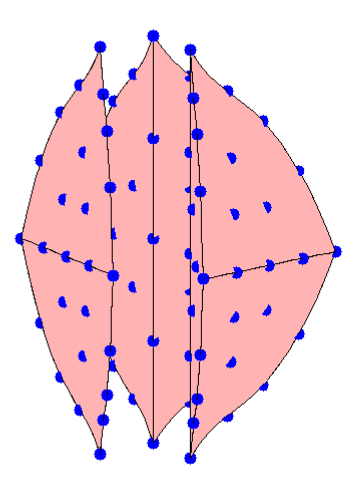} 
    \hspace{2mm}
    \includegraphics[width=0.23\textwidth]{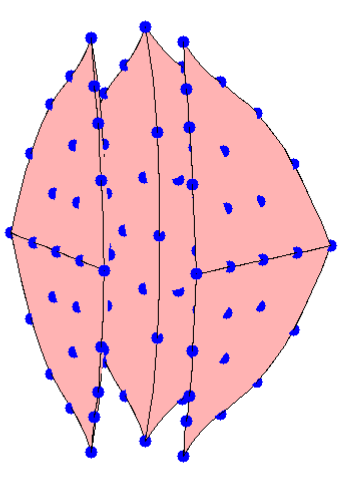} 
    \caption{Face swaps in 3D: original elements \textit{(left)}, initial straight-sided swap \textit{(middle)}, curved swap \textit{(right)}.}
    \label{fig:3dflip_schem}   
    \includegraphics[width=0.2\textwidth]{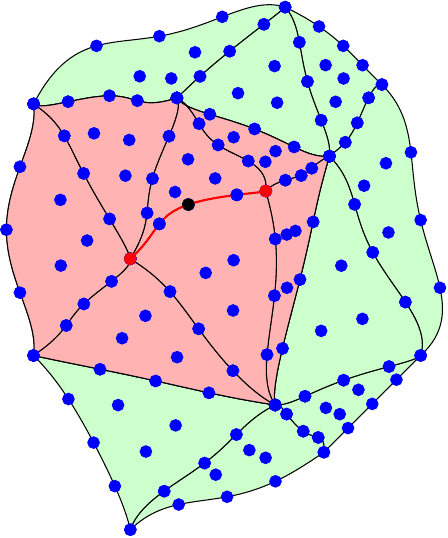} 
    \hspace{2mm}
    \includegraphics[width=0.2\textwidth]{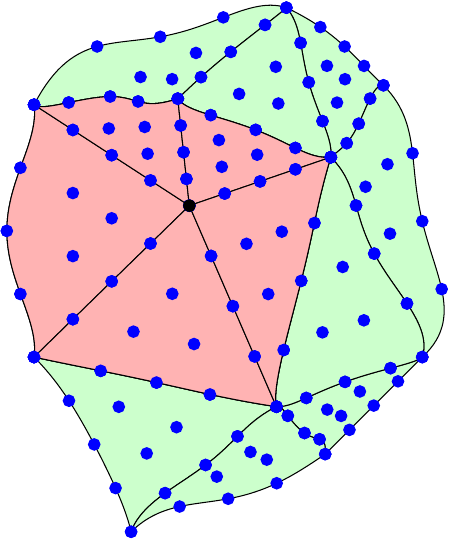} 
    \hspace{2mm}
    \includegraphics[width=0.2\textwidth]{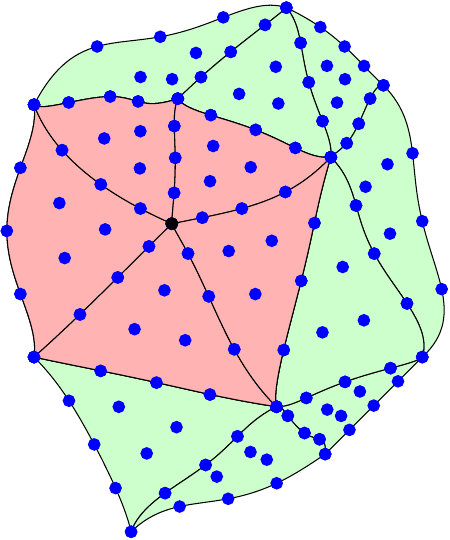} 
    \caption{Edge collapses in 2D: original elements \textit{(left)}, initial straight-sided collapse \textit{(middle)}, curved collapse \textit{(right)}.}
    \label{fig:2dclps_schem}    
    \includegraphics[width=0.25\textwidth]{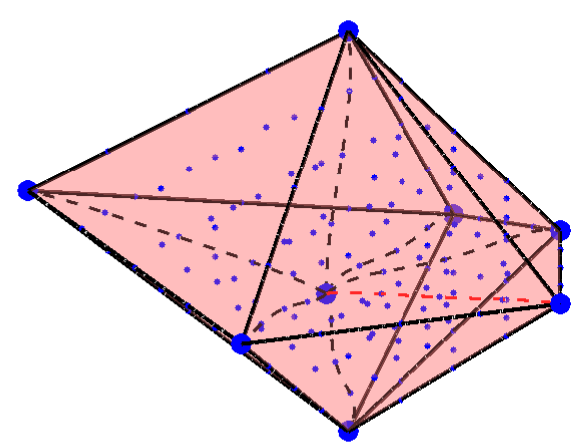} 
    \hspace{2mm}
    \includegraphics[width=0.25\textwidth]{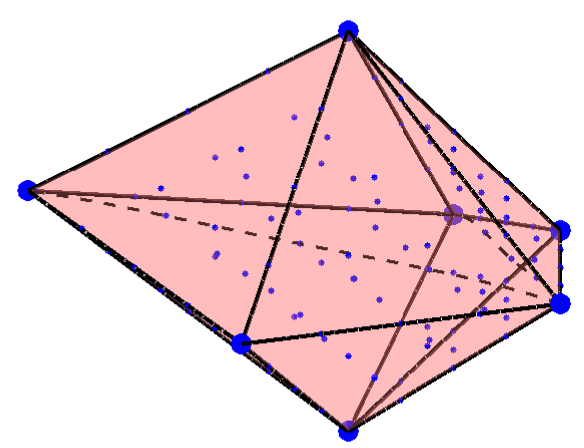} 
    \hspace{2mm}
    \includegraphics[width=0.25\textwidth]{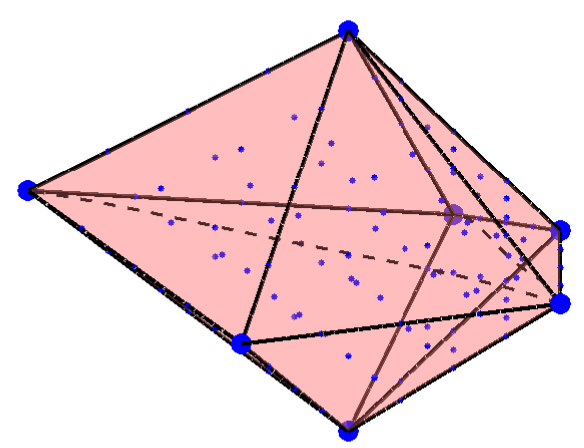} 
    \caption{Edge collapses in 3D: original elements \textit{(left)}, initial straight-sided collapse \textit{(middle)}, curved collapse \textit{(right)}.}
    \label{fig:3dclps_schem}        
\end{figure}

Once a group of elements is identified to perform edge/face swaps on,
we first perform the straight-sided initial guess.
This is done by freezing the boundary, 
discarding all the interior nodes and replacing them in the 
flipped configuration by the straight-sided information. 
This initial flip is extremely prone to tangling
even when starting out with fairly well-shaped elements,
which is more severe for the 3D case.
By the high-order distortion measure (\ref{EQ: HIGH ORDER DISTORTION}),
many high-order elements which seem valid by visual inspection
with all interior nodes contained within the boundary may very well be invalid.
Finally, we apply the simultaneous untangling and smoothing procedure
described in the previous section to arrive at the curved flip
(Figures \ref{fig:2dflip_schem} and \ref{fig:3dflip_schem}). 
We decide to accept the resulting flip if it results in a 
lower aggregate distortion $\hat{\eta}_{agg}$ from 
the original configuration.
The normalization by the number of elements in the group
in the definition is to allow comparison 
of the original and flipped configuration, since 
in 3D there are many topological changes
which in general may alter the number of elements.

\subsection{Edge Collapses}
Once an edge is identified to collapse,
we first do the straight-sided initial guess.
This is done by identifying the patch of elements that contains
either endpoint of the edge being collapsed, freezing its boundary, 
discarding all the interior nodes, 
and replacing them by the straight-sided edge collapse
to the midpoint of the original curved edge.
Then we apply the simultaneous untangling and smoothing procedure
to arrive the final curved collapse
(Figures \ref{fig:2dclps_schem} and \ref{fig:3dclps_schem}). 
If the edge identified for 
collapse is on the boundary of the domain, information will
inevitably be lost but can be mitigated in the high-order 
case by performing projection based interpolation
on the patch boundary by standard techniques encountered in $hp$-FEM methods \cite{demkowicz2006computing}.

\subsection{Edge Splitting}
\begin{figure}[hp]
    \centering
    \includegraphics[width=0.25\textwidth]{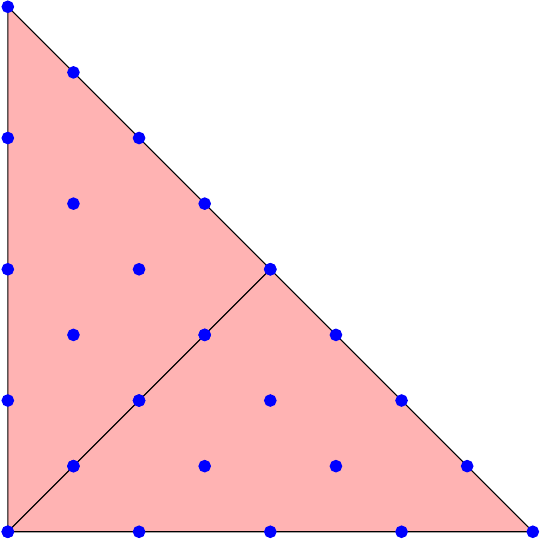} 
    \hspace{10mm}
    \includegraphics[width=0.25\textwidth]{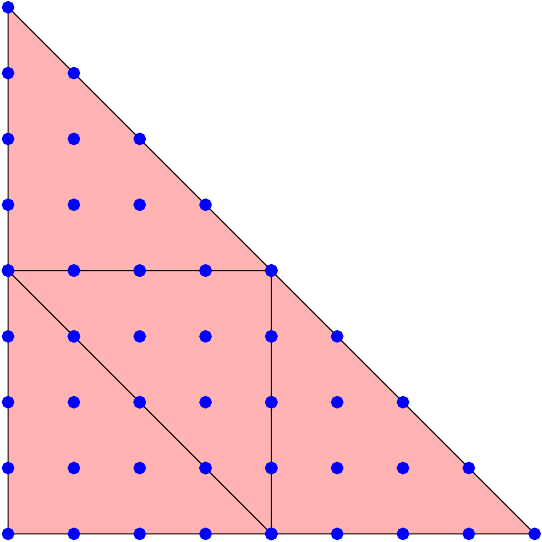} 
    \caption{Element subdivision templates in 2D: $1 \rightarrow 2$ \textit{(left)}
    and $1 \rightarrow 4$ \textit{(right)}.}
    \label{fig:2dsplit_templates}
    \includegraphics[width=0.3\textwidth]{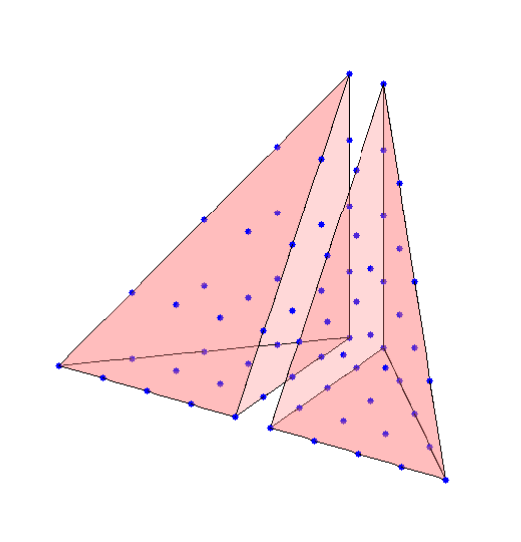} 
    \includegraphics[width=0.3\textwidth]{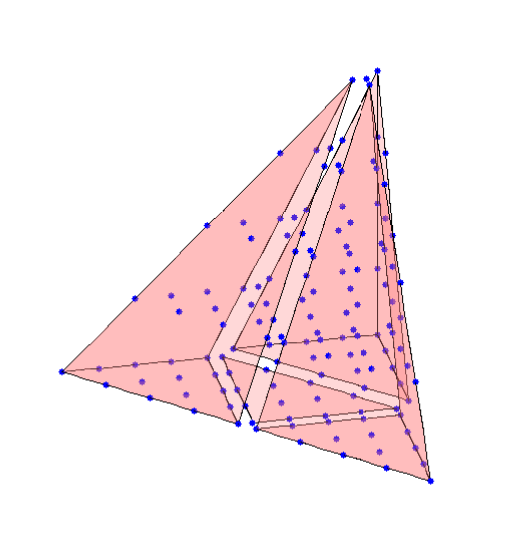}
    \includegraphics[width=0.3\textwidth]{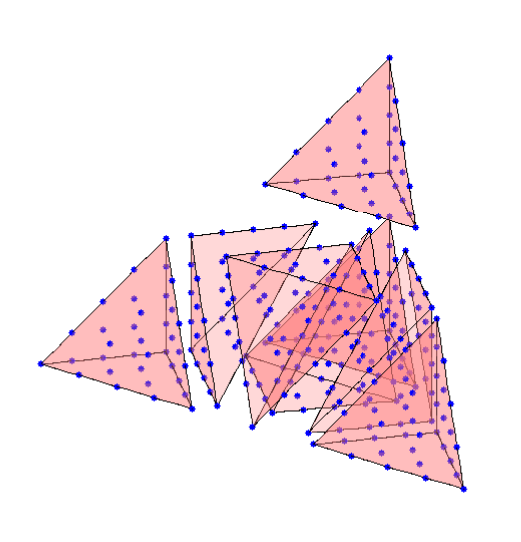}
    \caption{Element subdivision templates in 3D: $1 \rightarrow 2$ \textit{(left)}, $1 \rightarrow 4$ \textit{(middle)}, and $1 \rightarrow 8$ \textit{(right)}.}
    \label{fig:3dsplit_templates}
    \includegraphics[width=0.25\textwidth]{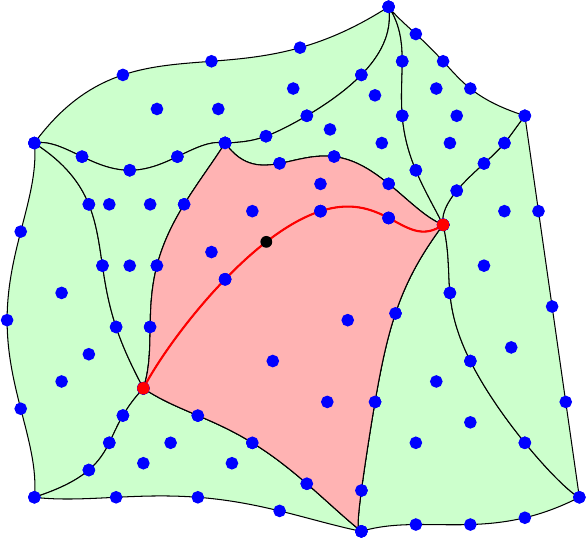} 
    \hspace{2mm}
    \includegraphics[width=0.25\textwidth]{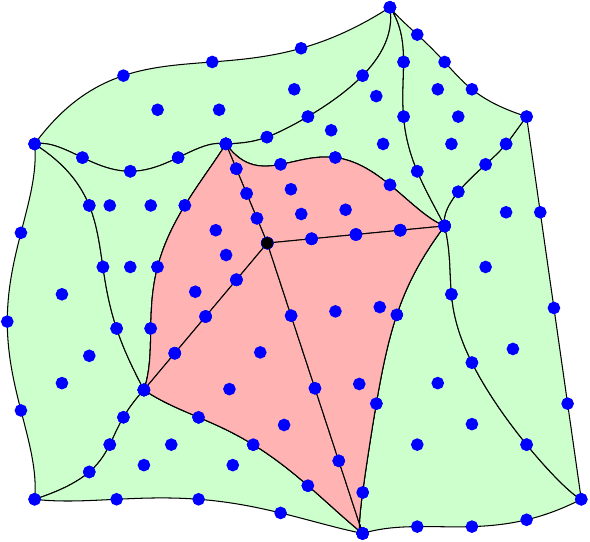} 
    \hspace{2mm}
    \includegraphics[width=0.25\textwidth]{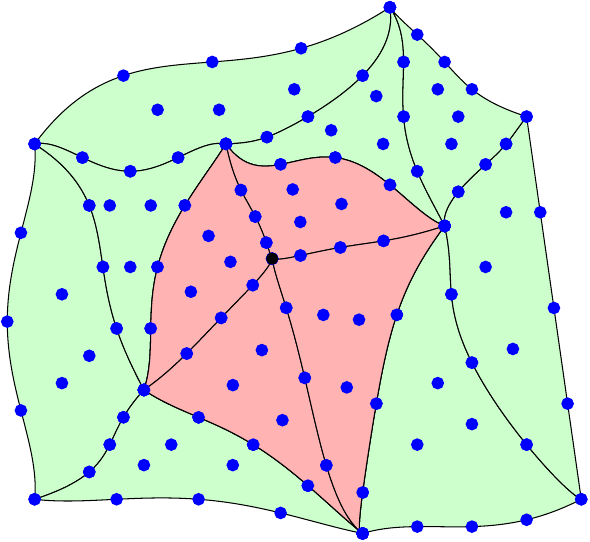} 
    \caption{Edge splitting in 2D: original elements \textit{(left)}, initial straight-sided split \textit{(middle)}, curved split \textit{(right)}.}
    \label{fig:2dsplit_schem}
    \centering
    \includegraphics[width=0.3\textwidth]{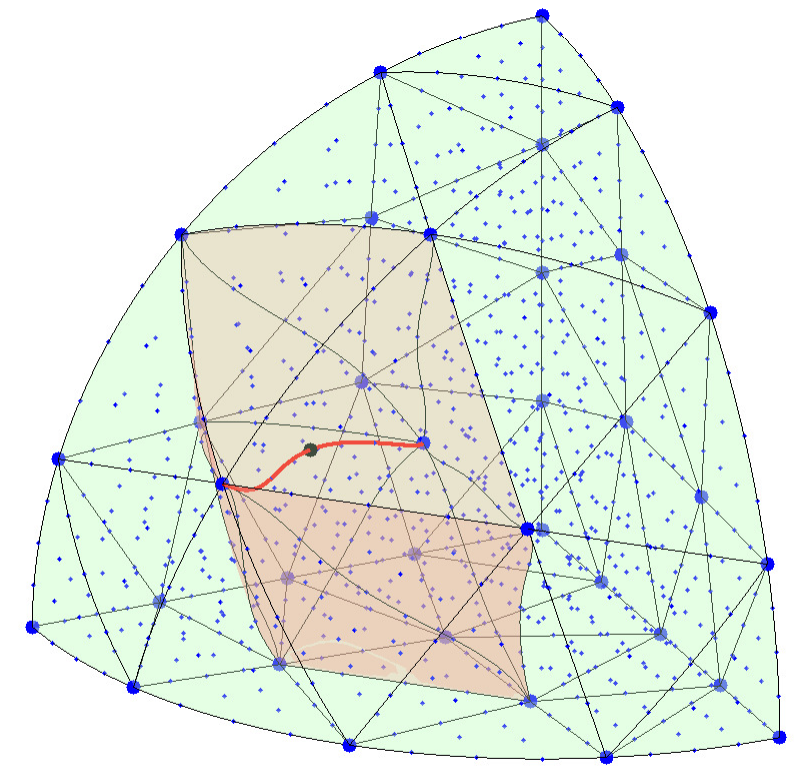} 
    \hspace{2mm}
    \includegraphics[width=0.3\textwidth]{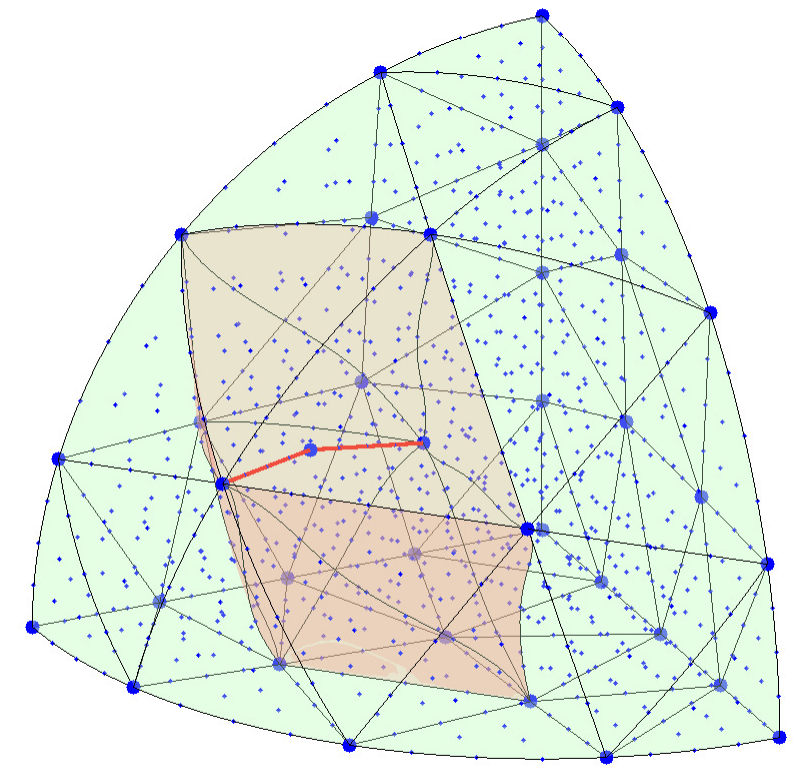} 
    \hspace{2mm}
    \includegraphics[width=0.3\textwidth]{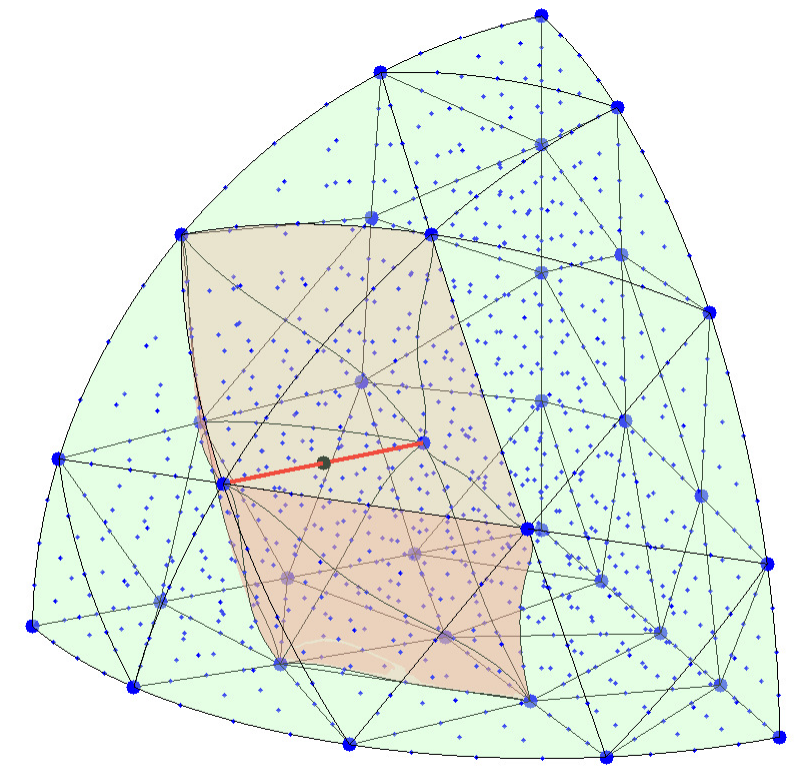} 
    \caption{Edge splitting in 3D: original elements \textit{(left)}, initial straight-sided split \textit{(middle)}, curved split \textit{(right)}.}
    \label{fig:3dsplit_schem}       
\end{figure}

Depending on the setting, the edges to be split could
either be directly identified or implied 
by choosing elements. 
Either way, once an element is identified 
to be split and the straight sided guess
is performed according to one of the templates 
in 2D (Figure \ref{fig:2dsplit_templates})
or 3D (Figure \ref{fig:3dsplit_templates}),
its neighbors must also be split by a straight-sided
guess with the corresponding template to prevent hanging nodes.
Then, we apply simultaneous untangling and smoothing
on this group of elements to arrive at the final curved split
(Figures \ref{fig:2dsplit_schem} and \ref{fig:3dsplit_schem}).
There are many more possible subdivision templates
in 3D depending on the number of marked edges 
(Figure 11 of \cite{de1999parallel}), and that for any of these templates
to determine how to split the edge(s), 
the amount of ``smoothing" required after the edge splits applied to the straight sided guess can differ quite greatly and alter the overall efficiency of the method, but this is outside
of the scope of this work.

%%%%%%%%%%NUMERICAL EXPERIMENTS
\section{NUMERICAL EXPERIMENTS} 
We demonstrate how local element operations for curved meshes can be combined 
with high-order mesh smoothing to maintain element quality and sizing 
for high-order meshes subject to severe deformations in both 2D and 3D.
All the meshes used for the numerical experiments are generated 
by Gmsh \cite{geuzaine2009gmsh} and are of polynomial degree four. \\

All the results have been obtained on a MacBook Pro with one dual-core Intel Core i7 CPU
(Apple  Inc.,  Cupertino,  California,  US), with a clock frequency of 3.0 GHz, and  total  memory of 16 GBytes.  
As a proof of concept, this code has been fully developed in MATLAB without using any additional toolbox. 
The code  is  not optimized, not  parallel, and not compiled and the results are provided as an illustration of this concept.

\subsection{Rotation}
We first demonstrate how node smoothing combined with edge flips
can maintain element quality under rotation. 
The 2D mesh is a circle-in-square mesh of 116 elements,
where the circle is centered at the origin with radius 0.5
inside the square $[-1, 1]^2$ (Figure \ref{fig:rot2d} \textit{left}). 
The analogous 3D mesh is a sphere-in-cube mesh of 3395 elements,
where the sphere is centered at the origin with radius 0.5
inside the cube $[-1, 1]^3$ (Figure \ref{fig:rot3d} \textit{left}). 
At each step, we rotate the nodes on the circle/sphere,
perform high-order node smoothing
while keeping the nodes  on the circle/sphere fixed, 
apply curved flips according to Algorithm \ref{alg:flipsweep}, 
and smooth again.
Algorithm \ref{alg:flipsweep} constructs a priority 
queue of elements with higher distortion than a
floating threshold value and considers all topologically possible
flips for each element, accepting the one (if any) that results
in the greatest improvement. It performs sweeps over 
all the elements with a reset floating threshold for each 
sweep until no further improvement is possible.

Rotation cannot exceed much more than a quarter turn with only node smoothing 
(Figures \ref{fig:rot2d} and \ref{fig:rot3d} \textit{middle}).
A straight sided mesh cannot rotate this far with smoothing alone; 
the presence of high order nodes grant us the flexibility to maintain valid 
(albeit highly stretched) elements for more severe deformations.
It is well known in the straight-sided case that 
such greedy flip algorithms like the one used here 
converge to an optimal triangulation (in the sense of Delaunay) 
in the 2D case, but offers no such guarantees in higher dimensions \cite{de2010structures}.
2D rotation can continue indefinitely,
and we can reasonably attain a good aggregate quality 
under rotation for the 3D case, but some poorly shaped elements 
emerge as evidenced by the growing maximum element distortion (Table \ref{table:rot}); 
this issue is no doubt magnified by our use of very coarse meshes.
We also report the aggregate mesh quality (\ref{EQ: AGG}), and the numerical
results also illustrate a fundamental issue of mesh smoothing discussed 
at the beginning of Section \ref{sec:optim}. It turns out even if the
aggregate distortion of the mesh is very low, this can obscure the fact that
there are elements of very low quality present as well. 
This can also be remedied with more sophisticated
mesh cleanup procedures featuring more specialized 
element operations such as compound operations \cite{joe1995construction}
and sliver exudation \cite{cheng2000silver}, 
but such issues are not the focus of this work. 
Also note that here we do not nodes on the boundary or the circle/sphere. 
Allowing these nodes to slide on the fixed boundaries while smoothing
the entire mesh would certainly be an improvement to the procedure.

\begin{table}[!ht]
    \centering 
    
    \includegraphics[width=0.3\textwidth]{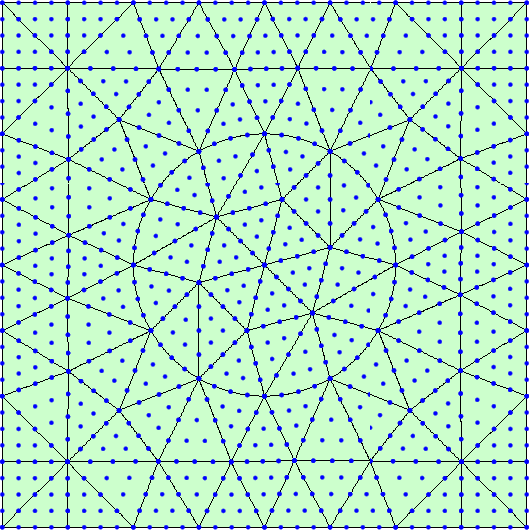}
    \hspace*{0.5cm}
    \includegraphics[width=0.3\textwidth]{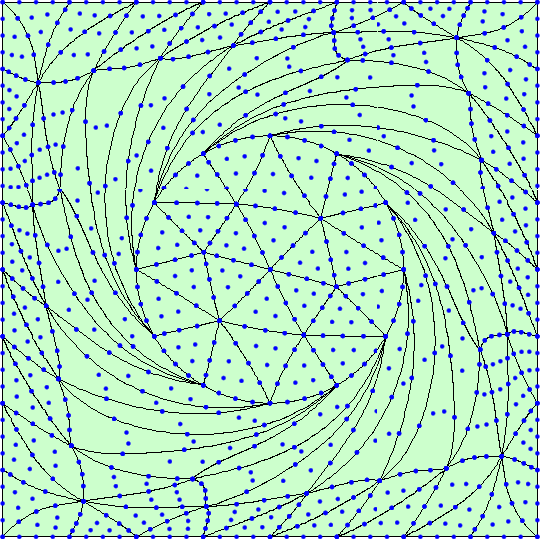} 
    \hspace*{0.5cm}
    \includegraphics[width=0.3\textwidth]{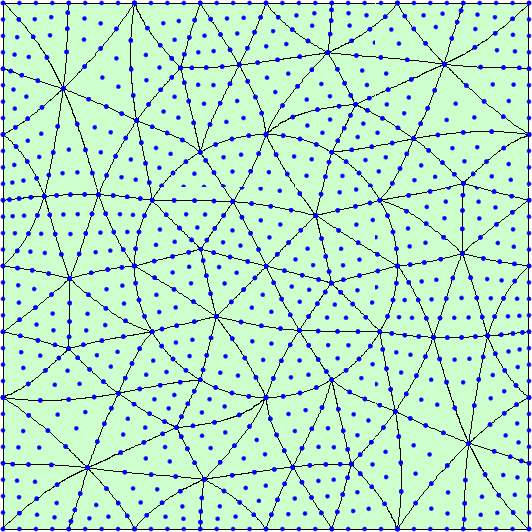}   
    \captionof{figure}{2D rotation with flips: original mesh (\textit{left}), quarter turn with smoothing only (\textit{middle}), quarter turn with smoothing and flips (\textit{right}).} 
    \label{fig:rot2d}
    \vspace{1cm}
    
    \includegraphics[width=0.3\textwidth]{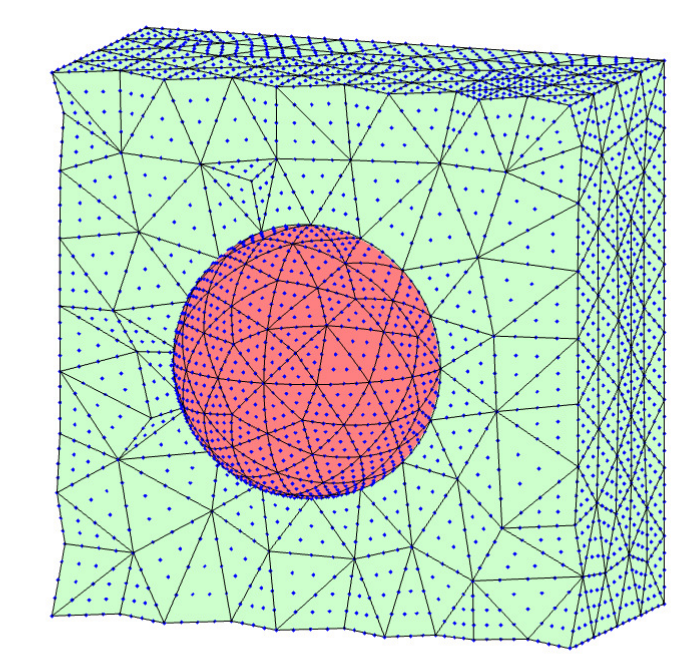}
    \hspace*{0.5cm}
    \includegraphics[width=0.3\textwidth]{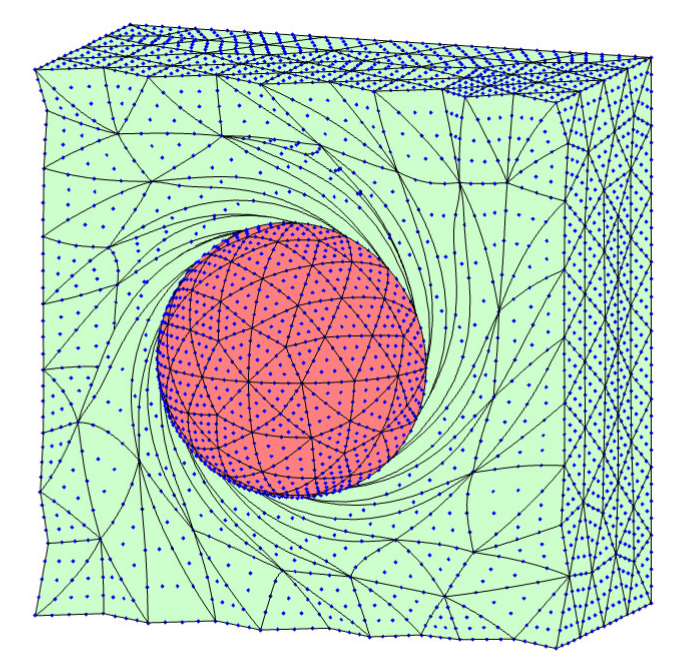} 
    \hspace*{0.5cm}
    \includegraphics[width=0.3\textwidth]{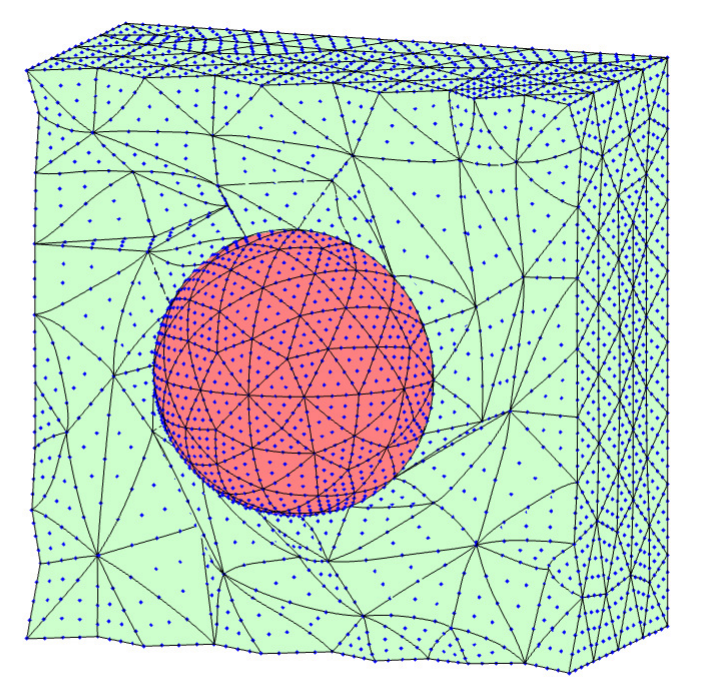}   
    \captionof{figure}{3D rotation with flips: original mesh (\textit{left}), quarter turn with smoothing only (\textit{middle}), quarter turn with smoothing and flips (\textit{right}).}
    \label{fig:rot3d} 
    
    \captionof{table}{Mesh Distortion of Figures \ref{fig:rot2d} and \ref{fig:rot3d}}   
    \label{table:rot} 
    \begin{tabular}{ccccc}
    & \multicolumn{2}{c}{2D} & \multicolumn{2}{c}{3D} \\
    \cline{2-5}
    & max & agg & max & agg \\
    \hline
    Original mesh  & 1.167 & 1.141 & 2.499 & 1.669 \\
    Smoothing only & 7.876  & 5.327 & 15.531 & 8.490 \\
    With flips  & 1.197 & 1.179 & 4.086 & 1.826 \\
    \hline
    \end{tabular}  
\end{table}

\subsection{Coarsening}

\begin{figure}
 \centering
  \includegraphics[width=0.4\textwidth]{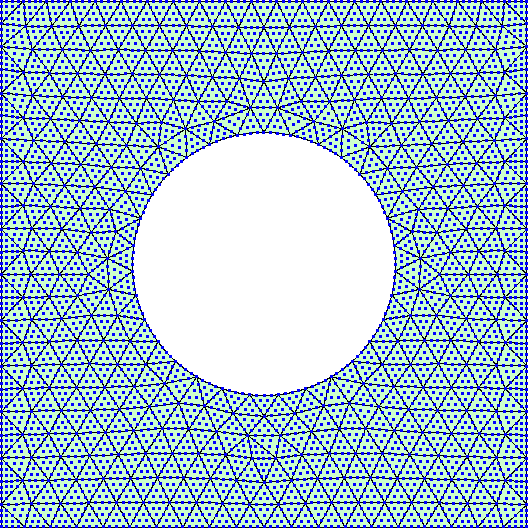}
  \hspace*{0.5cm}
  \includegraphics[width=0.4\textwidth]{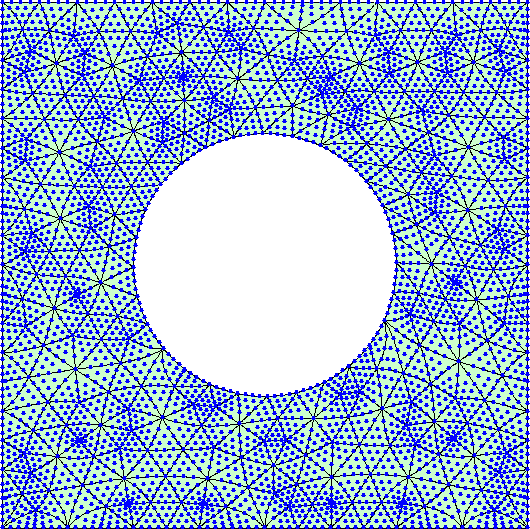} \\
  \vspace*{0.5cm}
  \includegraphics[width=0.4\textwidth]{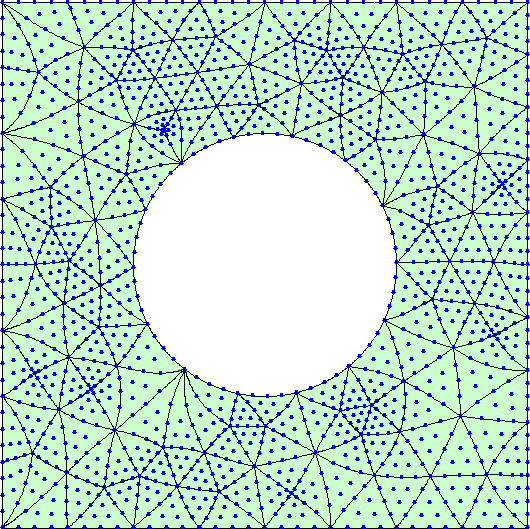}   
  \hspace*{0.5cm}
  \includegraphics[width=0.4\textwidth]{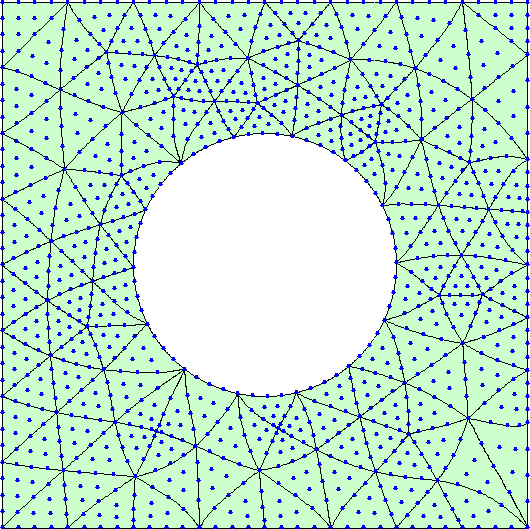}   
 \caption{2D coarsening with collapses: original mesh (\textit{top left}), one sweep of collapses (\textit{top right}), five sweeps of collapses (\textit{bottom left}), final mesh (\textit{bottom right}).}
 \label{fig:clps2d}
\end{figure}

Though this example does not deal with mesh deformation,
we demonstrate how curved collapses can be combined
with high-order node smoothing to coarsen a mesh (Figure \ref{fig:clps2d}).
We begin with a rather fine 2D mesh with 690 elements
of the square $[-1, 1] \times [-1, 1]$ with a 
circular hole at the origin with radius 0.5 (Figure \ref{fig:clps2d}, \textit{top left}).
We apply multiple passes of curved edge collapses to coarsen the mesh
until convergence according to a chosen uniform size function 
(Algorithm \ref{alg:clpssweep}). In this case we choose the ideal edge length
$L_0 = 0.1$ as a parameter far enough from the original edge length just to 
illustrate the procedure to induce collapses.
We also coarsen elements on the interior boundary 
with projection by interpolation in the $L^2$
error minimizing sense to maintain an accurate 
representation of the circular boundary with fewer elements.

\subsection{Translation}

We demonstrate the combination of all three
curved element operations with mesh smoothing 
to maintain both element quality and uniform sizing. 
The 2D mesh is a circle-in-rectangle mesh of 314 elements,
where the circle is centered at $(-1, 0)$ and has radius 0.5 
inside the rectangle $[-2, 2] \times [-1, 1]$. 
The analogous 3D mesh is a sphere-in-prism mesh of 3395 elements,
where the sphere is centered at $(-1, 0, 0)$ with radius 0.5
inside the prism $[-2, 2] \times [-1, 1] \times [-1, 1]$. 
At each step, we translate the nodes on the circle/sphere,
perform high-order node smoothing while keeping the nodes 
on the circle/sphere fixed, apply all three types
of local mesh operations according to 
Algorithm \ref{alg:allsweep}, and smooth again. 
Algorithm \ref{alg:allsweep} applies sweeps of splits (Algorithm \ref{alg:splitsweep}) and collapses (Algorithm \ref{alg:clpssweep}) until convergence, 
and then applies sweeps of flips (Algorithm \ref{alg:flipsweep}) until convergence.
In terms of the parameters for both Algorithms \ref{alg:clpssweep} and \ref{alg:splitsweep},
the ideal edge length is set to be $L_0 = \text{avg}(L)$,
where $L$ is the list of computed original edge lengths.
This choice is made simply because after the translation, we want 
the final mesh to resemble our original mesh, which 
was roughly uniform and with the same element sizes. As an aside,
one could imagine introducing anisotropy to the mesh
through these local element operations through a specified size function
which would specify the ideal lengths $L_0$ at any point in space, but
this is beyond the scope of the work.

\begin{figure}[htbp!]
    \centering 
    \includegraphics[width=0.48\textwidth]{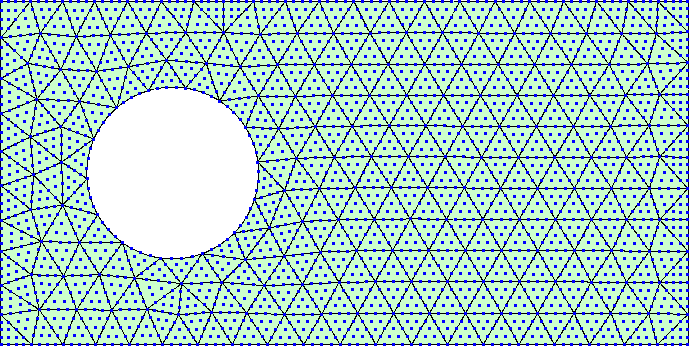}
    \hspace*{0.25cm}
    \includegraphics[width=0.48\textwidth]{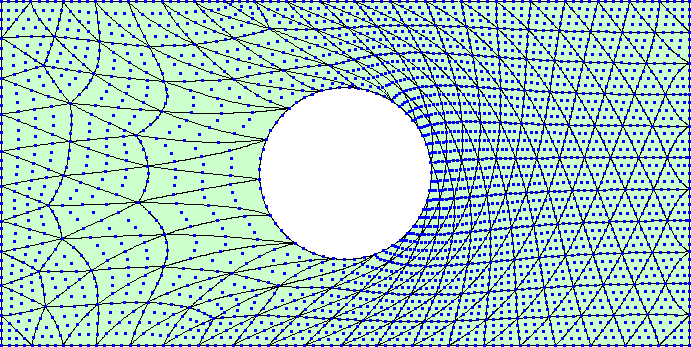} \\
    \vspace*{0.25cm}
    \includegraphics[width=0.48\textwidth]{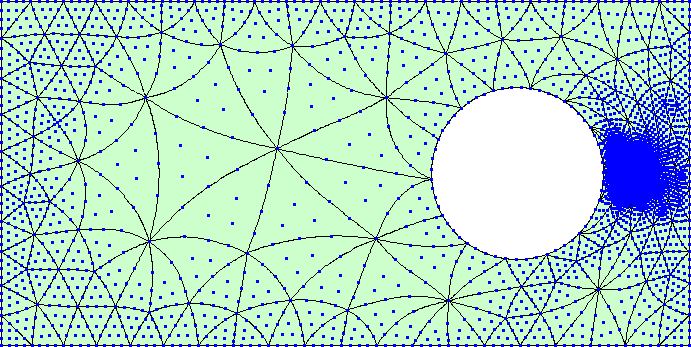}   
    \hspace*{0.25cm}
    \includegraphics[width=0.48\textwidth]{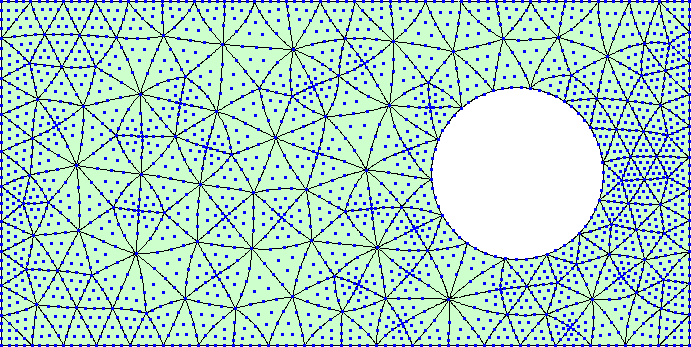}   
    \caption{2D translation: original mesh (\textit{top left}), smoothing only (\textit{top right}), smoothing with flips (\textit{bottom left}), smoothing with flips, collapses, and splitting (\textit{bottom right}).}
    \label{fig:trans2d}
    \vspace{0.5cm}
\end{figure}

\begin{figure}[htbp!]
    \includegraphics[width=0.48\textwidth]{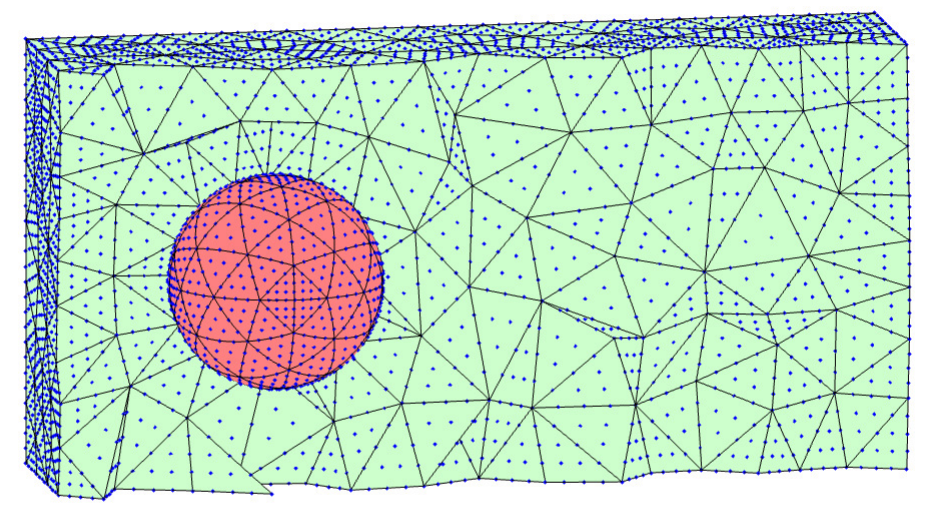}
    \hspace*{0.25cm}
    \includegraphics[width=0.48\textwidth]{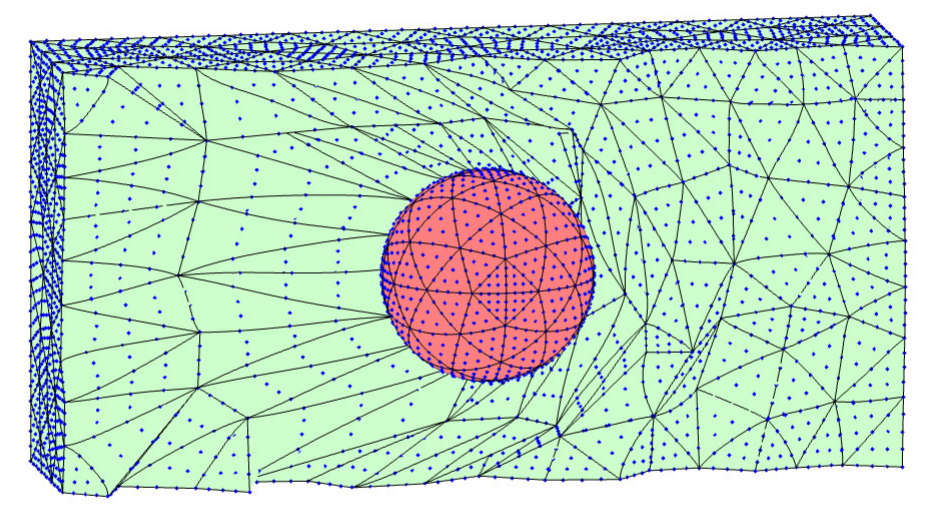} \\
    \vspace*{0.25cm}
    \includegraphics[width=0.48\textwidth]{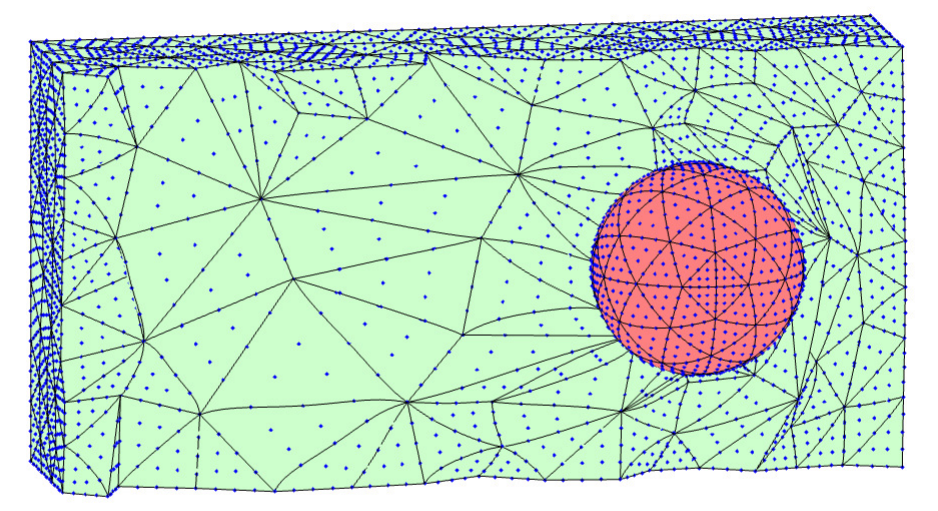}   
    \hspace*{0.25cm}
    \includegraphics[width=0.48\textwidth]{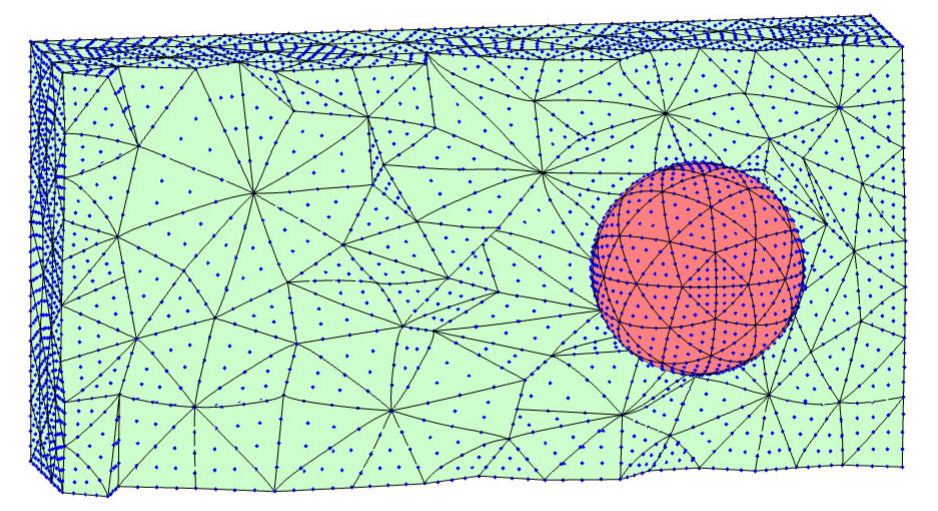}   
    \captionof{figure}{3D translation: original mesh (\textit{top left}), smoothing only (\textit{top right}), smoothing with flips (\textit{bottom left}), smoothing with flips, collapses, and splitting (\textit{bottom right}).}
    \label{fig:trans3d} 
    \centering     
    \captionof{table}{Mesh Distortion of Figures \ref{fig:trans2d} and \ref{fig:trans3d}}
    \label{table:trans}
    \begin{tabular}{ccccc}
    & \multicolumn{2}{c}{2D} & \multicolumn{2}{c}{3D} \\
    \cline{2-5}
     & max & agg & max & agg \\
    \hline
    Original mesh  & 1.1821 & 1.0473 & 2.2466 & 1.7481 \\
    Smoothing only  & 2.8378 & 2.3771 & 7.3829 & 3.8409 \\
    With flips  & 2.0643 & 1.2075 & 3.8427 & 2.1742 \\
    Flips, collapses, and splitting & 1.1821 & 1.0498 & 3.0003 & 1.9666 \\
    \hline
    \end{tabular}    
\end{figure}

Translation cannot make it more than halfway across the 
prism with only node smoothing (Figures \ref{fig:trans2d} and \ref{fig:trans3d}).
With the addition of flips, we can make it across the entire domain,
but this results in elements of highly variable size.
This can be remedied with the addition of collapses 
and splitting to maintain a uniform size, 
which allows us to practically maintain the original 
element sizes and qualities after
translation across the entire domain (Table \ref{table:trans}).

\section{CONCLUSION AND FUTURE WORK}
We extend the local element operations of edge flips, 
edge collapses and edge splitting to the case of curved meshes.
Our numerical examples demonstrate how these local operations can be
combined with node smoothing to successfully
perform mesh adaptation directly on curved meshes to maintain
high element quality and desired sizing. 
Curved meshes can naturally deal with more severe deformations 
than straight sided ones on much coarser meshes
due to the higher number of degrees of freedom, and
mesh adaptation can be applied to take advantage 
of this fact to deal with more severe deformations. 
These operations are more expensive to perform in the curved setting;
for straight-sided meshes they amount to a simple permutation 
of the existing data and require no node smoothing.
However, in the big picture this is a small price to pay,
since in any practical setting these local operations
deal with a relatively small number of elements 
when used as an alternative to a 
far more costly procedure such as remeshing.

We see two avenues of future work to fully unlock the utility
of these operations towards mesh adaptation for curved meshes.
More efficient mesh optimization procedures combining these
ingredients need to be developed; in this work, 
we use fairly naive and greedy algorithms which are
straightforward, but far from optimal.
As noted in \cite{geuzaine2009gmsh}, mesh optimization procedures
currently are based on ``black magic: even if the ingredients 
required to construct a mesh optimization procedure are well known ...
there is no known best recipe, i.e., no known optimal way of combining these." 
In the 15 years since this statement was made, there have
certainly been many improvements in the ``ingredients", 
but there doesn't seem to have much work focused solely on the ``recipes". 
The large body of existing work on these mesh 
optimization procedures is mainly concerned with achieving the
best possible mesh quality without much concern for speed.
This is a significant limitation in the curved setting,
since node smoothing is the main bottleneck in these procedures
and is much more expensive in the high-order case.
As Klinger and Shewchuk note in \cite{klingner2008aggressive},
mesh \textit{optimization} procedures have been shown 
to attain results as good as those from any mesh \textit{generation} procedure,
and ``if the barrier of speed can be overcome, 
the need to write separate programs for mesh generation
and mesh improvement might someday disappear."
This is particularly relevant for the issue of high-order
mesh generation, since the limited number of procedures 
currently in use are nearly all of the \textit{a priori}
variety rather than direct methods. Surely, this is in 
no small part due to our limited ability to
currently perform mesh adaptation directly on curved meshes.

The real challenge in the high-order setting is the
ability to perform smoothing quickly and untangling robustly. 
It is worth emphasizing just how much of the cost 
comes from the smoothing procedure.
In this work we prescribe the mesh motion,
which must be followed by a smoothing of the entire mesh,
a pass of local mesh operations, and 
and then another smoothing of the entire mesh. 
The two contributions to the cost of each local mesh operation
is a simple permutation of indices in the mesh data structure for the straight sided guess
followed by smoothing. Since the cost of the former is practically negligible,
virtually all of the cost comes from the mesh smoothing procedure,
made even more expensive in a \textit{high-order} mesh. 
While the approach for untangling considered in this work of the
regularization of a mesh distortion metric
has the advantage that it is relatively simple, 
it is very dependent on a heuristically chosen parameter
and many authors have pointed out that it
fails quite often in practice for various reasons, 
such as if the boundary of the mesh is also tangled or if
the inversion primarily occurs at the corners of the elements
\cite{gargallo2015optimization, irving2004invertible}. 
Another limitation of these optimization-based 
approaches is that they are largely based on the special choice
or modification of some objective function and lack 
specialized solvers. A promising idea is the extension
of some of the ideas in \cite{fortunato15winslow, knupp1999winslow, karman2006mesh} towards a fully unstructured, 
high-order Winslow-based smoothing and untangling procedure. 
Winslow-based methods have been shown to be 
very robust in practice and do not require the choice 
of any empirical parameters. Each iteration of these methods
essentially requires an elliptic solve, for which there 
exists highly robust and efficient solvers (i.e. multigrid)
and would make a huge difference in efficiency for the high-order setting.

From here, we can revisit the vast amount of literature
on mesh optimization procedures almost exclusively for the 
straight-sided case to better understand the differences
and limitations in the curved setting, 
which we know very little about.
Further developments in this area would also remove
one of the constraints on the use of
high-order elements for the so called 
arbitrary Lagrangian-Eulerian (ALE) simulation framework 
for moving meshes with large deformations \cite{wang2015discontinuous}. \\

\textbf{Data Availability Statement:} Research data are not shared.

\newpage
\appendix
\begin{algorithm}[ht]
\DontPrintSemicolon
\SetAlgoLined
\lIf{dim = 2D} {$tmin = 1.5$} \lElseIf{dim = 3D} {$tmin = 2$} 
Input curved mesh $M$ with elements $\{e_i\}$ \\
Set distortion threshold $T = \min(tmin, \min(D)), D(i) = $ distortion of $e_i$ \;
Create priority queue $P = \{(e_j, D(j)), \forall j \ D(j) > T \text{ descending} \}$\\
 \While{true}{
  \For{i in P}{
  \lIf{dim = 2D} {consider all possible $2 \rightarrow 2$ flips with $e_i$} 
  \lElseIf{dim = 3D} {consider all possible $2 \rightarrow 3, 3 \rightarrow 2, 4 \rightarrow 4$ flips with $e_i$}  
  \If{best of topologically possible flips improves $\hat{\eta}_{agg}$}{
   Update the mesh $M$ \\
   Remove all involved elements from the queue $P$ and renumber
   }
  }  
  \If{no flips done in sweep}{
   \Return $M$\;
   }
  Recompute $\{e_i\}$, $T$ and $P$
 }
 \caption{Application of Curved Flips}
 \label{alg:flipsweep} 
\end{algorithm}

\begin{algorithm}[ht]
\DontPrintSemicolon
\SetAlgoLined
Set ideal edge length as parameter $L_0$ \\
Input curved mesh $M$ with edges $\{E_j\}$ \\
Compute list of edge lengths $L$ \\
Create priority queue $P = \{(E_k, L(E_k)), \forall k \ \frac{L_0}{L(E_k)} > 1.5 \text{ descending} \}$\\
 \While{true}{
  \For{i in P}{
  Perform curved collapse on edge $e_i$ and update the mesh $M$ \\
  Remove all edges in collapsed patch from $P$ and renumber
   }
  \If{no collapses done in sweep}{
   \Return $M$\;
   }
  Recompute $\{E_j\}$, $L$ and $P$
 }
 \caption{Application of Collapses}
 \label{alg:clpssweep} 
\end{algorithm}

\begin{algorithm}[ht]
\DontPrintSemicolon
\SetAlgoLined
Set ideal edge length as parameter $L_0$ \\
Input curved mesh $M$ with edges $\{E_j\}$ \\
Compute list of edge lengths $L$ \\
Create priority queue $P = \{(E_k, L(E_k)), \forall k \ \frac{L_0}{L(E_k)} > 1.5 \text{ descending} \}$\\
 \While{true}{
  \For{i in P}{
  Identify elements $K$ containing edge $E_i$ \\
    \For{j in K}{
        Identify all edges of element $e_j$ in queue $P$ to choose template \\
        Perform curved split on element $e_j$ and update the mesh $M$ \\
        Remove all edges in collapsed patch from $P$ and renumber
    }
   }
  \If{no splits done in sweep}{
   \Return $M$\;
   }
  Recompute $\{E_j\}$, $L$ and $P$
 }
 \caption{Application of Splitting}
 \label{alg:splitsweep} 
\end{algorithm}

\begin{algorithm}[ht]
\DontPrintSemicolon
\SetAlgoLined
Input curved mesh $M$ with elements $\{e_i\}$ and edges $\{E_j\}$ \\
    Perform sweeps of splits (Algorithm \ref{alg:splitsweep})  \\
    Perform sweeps of collapses (Algorithm \ref{alg:clpssweep}) \\
    Perform smoothing of the entire mesh $M$ \\
    Perform sweeps of flips (Algorithm \ref{alg:flipsweep}) \\
    Perform smoothing of the entire mesh $M$ \\
    \Return $M$\;
 \caption{Application of all operations: flips, coarsening, and splitting}
 \label{alg:allsweep} 
\end{algorithm}

\newpage

\bibliographystyle{wileyj}
\bibliography{pap}

\end{document}